\centerline{\bf Real Regulators on Milnor Complexes}
\bigskip
\centerline{\tt (Some Informal Notes; December 2000)}
\bigskip
\centerline{James D. Lewis}
\bigskip
\S {\bf 0. Introduction} 
\bigskip
Let $X/_{\bf C}$ be a projective algebraic manifold
of dimension $n$, with corresponding sheaf of regular functions
${\cal O}_X$. Put 

$$
{\cal K}^M_{k,X} :=
\big(\underbrace{{\cal O}_X^{\times}\otimes \cdots\otimes
{\cal O}_X^{\times}}_{k-{\rm times}}\big) / J\quad\quad
{\rm (Milnor\ sheaf)},
$$
where $J$ is the subsheaf of the tensor product generated
by sections of the form $\tau_1\otimes \cdots \otimes
\tau_k$
such that $\tau_i + \tau_j = 1$ for some $i$ and $j$. 
Put  $\overline{\cal K}^M_{k,X} := \ {\rm Image}\big(j: {\cal K}^M_{k,X}
\to {K}_{k}^M({\bf C}(X))\big)$.
Our goal is to give an simple and explicit description of a regulator 
map to a certain quotient of real Deligne cohomology, 
in terms of logarithms (for $m\geq 1$):
$$
r_{\log} : H_{\rm Zar}^{k-m}(X,\overline{\cal K}^M_{k,X})
\to H_{\cal D}^{2k-m}(X,{\bf R}(k))/(?),
$$
where we recall that
$$ 
H_{\cal D}^{2k-m}(X,{\bf R}(k))\simeq {H^{2k-m-1}(X,{\bf R}(k-1))\over
\pi_{k-1}(F^kH^{2k-m-1}(X,{\bf C}))}.
$$
One of the reasons for constructing such a map, is the
relationship between $H_{\rm Zar}^{k-m}(X,{\cal K}^M_{k,X})$ and
$CH^k(X,m)$ (Bloch's higher Chow group [Blo]), and the
existence of a regulator
$r_{\cal D} : CH^k(X,m) \to H_{\cal D}^{2k-m}(X,{\bf R}(k))$.
Further, information about the regulator $r_{\log}$, 
or a suitable (Beilinson type) variant
$r_B : H_{\rm Zar}^{k-m}(X,{\cal K}^M_{k,X}) \to 
H_{\cal D}^{2k-m}(X,{\bf R}(k))/(?)$ should yield some information about
$r_{\cal D}$.
\bigskip
One has an exact sequence of sheaves ((Gabber, 1992), see [MS2]):
$$
{\cal K}^M_{k,X} \to {K}_{k}^M({\bf C}(X)) \to
\bigoplus_{{\rm codim}_XZ = 1}{K}_{k-1}^M({\bf C}(Z))\to\cdots
\to \bigoplus_{{\rm codim}_XZ = k-m}{K}_{m}^M({\bf C}(Z))\to\cdots
$$
$$
\to\bigoplus_{{\rm codim}_XZ = k-2}{K}_{2}^M({\bf C}(Z))
\to\bigoplus_{{\rm codim}_XZ = k-1}{K}_{1}^M({\bf C}(Z))\to
\bigoplus_{{\rm codim}_XZ = k}{K}_{0}^M({\bf C}(Z))\to 0,
$$
which clearly defines a flasque resolution of  $\overline{\cal K}^M_{k,X}$.
Similarly, for a quasi-projective variety $W$, 
if we set ${\cal CH}^k(r) =$ sheaf associated to the presheaf
$U\subset W$ open $\mapsto CH^k(U,r)$, then there is proven
in [Blo] a Gersten resolution:
$$
0\to {\cal CH}^k(r)\to \bigoplus_{x\in W^0}i_xCH^k({\rm Sp}(k(x)),r)\to
\bigoplus_{x\in W^1}i_xCH^{k-1}({\rm Sp}(k(x)),r-1)\to\cdots
$$
$$
\to \bigoplus_{x\in W^r}i_xCH^{k-r}({\rm Sp}(k(x)),0)\to 0.
$$
When $r=k$ (and $X=W$), both resolutions in fact agree (see [MS2]), and
thus one has ${\cal CH}^k(k) \simeq \overline{\cal K}^M_{k,X}$. Furthermore,
Bloch (op. cit.) constructs a spectral sequence,
$E_2^{p,q} := H^p_{\rm Zar}(X,{\cal CH}^k(-q))
\Rightarrow CH^k(X,-p-q)$. We now set $p+q=-m$, or $-q = p+m$.
Note that $p\geq 0,\ q\leq 0$. So $E_2^{p,-p-m} =
H^p(X,{\cal CH}^k(p+m))\Rightarrow CH^k(X,m)$.  
Now if we set $k=p+m$, then
$E_2^{k-m,-k} = H^{k-m}_{\rm Zar}(X,\overline{\cal K}^M_{k,X})$.
If $m\leq 2$ then one can argue by partial degeneration of
this spectral sequence that $CH^k(X,m) \simeq H_{\rm Zar}^{k-m}(X,
\overline{\cal K}^M_{k,X})$ ([MS] (op. cit.)).

\bigskip
We prove the following:
\bigskip
{\bf Main Theorem.} (i) The current defined by
$$
\left(\matrix{f_1,\ldots,f_m\in {\bf C}(Z)^{\times}\cr
{\rm codim}_XZ = k-m\cr}\right)
\ \mapsto \ \biggl(w \mapsto \sum_{\ell=1}^m\int_Z(-1)^{\ell-1}\log |f_\ell|
(d\log |f_1|\wedge\cdots \wedge \widehat{d\log |f_\ell|}\wedge\cdots 
\wedge d\log |f_m|)\wedge w\biggr),
$$
descends to a cohomological map 
$$
r_{\log} : H_{\rm Zar}^{k-m}
(X,\overline{\cal K}^M_{k,X})\to H_{\cal D}^{2k-m}(X,{\bf R}(k))/(?) := 
\big\{H^{k-1,k-m}(X)\oplus H^{k-m,k-1}(X)\big\}
\cap H^{2k-m-1}(X,{\bf R}(k-1)).
$$

(ii) There is a similar (but more complicated, see (3.0) below) 
explicit description for the composite 
$$
H_{\rm Zar}^{k-m}(X,\overline{\cal K}^M_{k,X}) \ {\buildrel r_B \over \to}
\big\{H^{k-1,k-m}(X)\oplus H^{k-m,k-1}(X)\big\}
\cap H^{2k-m-1}(X,{\bf R}(k-1)).
$$

(iii) In the case $m \leq 2$, both maps define regulators
$r_{\log} : CH^k(X,m) \to H_{\cal D}^{2k-m}(X,{\bf R}(k))$,
$r_B = r_{\cal D} : CH^k(X,m) \to H_{\cal D}^{2k-m}(X,{\bf R}(k))$ 
which agree for $m = 1$, and up to composite (for $m=2$)
with the real operator  
$$
J \in {\rm End}_{\bf R}\big(\big\{H^{k-1,k-m}(X)\oplus H^{k-m,k-1}(X)\big\}
\bigcap H^{2k-m-1}(X,{\bf R}(k-1))\big),
$$
$$ w= w_{k-1,k-m}\oplus w_{k-m,k-1}\quad \mapsto \quad
\sqrt{-1}w_{k-1,k-m}\ominus \sqrt{-1}w_{k-m,k-1}.
$$
\vskip.25in
\S {\bf 1. Review of Milnor $K$-Theory}
\bigskip
We first recall the definition
of Milnor $K$-theory [B-T]. Let ${\bf F}$ be a field
with multiplicative group ${\bf F}^{\times}$, and set:
$$
T({\bf F}^{\times}) = \coprod_{n\geq 0}T^n({\bf F}^{\times}),
$$
the tensor product of the ${\bf Z}$-module ${\bf F}^{\times}$. Thus
${\bf F}^{\times} \ {\buildrel\sim\over \to}\ 
T^1({\bf F}^{\times})$, $a\mapsto [a]$.
If $a \ne 0,\ 1$, set $r_a = [a]\otimes [1-a] \in T^2({\bf F}^{\times})$.
The two-sided ideal $R$ generated by $r_a$ is graded, and we put:
$$
K_{\ast}^M{\bf F} = T({\bf F}^{\bullet})/R = \coprod_{n\geq 0}K^M_n{\bf F}.
$$
Thus $K_{\ast}^M{\bf F}$ is presented as a ring by generators
$\ell(a)$ ($a\in {\bf F}^{\times}$) subject to:
\bigskip
$
(R_1)\hskip1in \ell(ab) = \ell(a) + \ell(b), 
$
\bigskip
$
(R_2) \hskip1in \ell(a)\ell(b) = 0\ {\rm if}\ a+b=1,\quad (a\ne 0,1).
$
\bigskip 
We summarize those results we need from [B-T].
\bigskip
{\bf Proposition 1.0.} The following relations are a consequence
of $R_1$ and $R_2$ above:
\bigskip
$
\hskip1in (R_3)\hskip1in \ell(a)\ell(-a) = 0.
$
\bigskip
$
\hskip1in (R_3^{\prime})\hskip1in \ell (a)\ell(-1) = -(\ell(a))^2.
$
\bigskip
$
\hskip1in (R_4) \hskip1in \ell(a)\ell(b) = -\ell(b)\ell(a).
$\quad
[Thus $K^M_{\ast}{\bf F}$ is anti-commutative.]
\bigskip
$
\hskip1in (R_5)\hskip1in \ell(a_1)\cdots \ell(a_n) = 0,\ {\rm if}\
a_1+\cdots+a_n = 1\ {\rm or}\ 0, \ {\rm and}\ n\geq 2.
$
\bigskip
Furthermore, ${\bf Z} \ {\buildrel \sim \over\to}\
K^M_0({\bf F})$, $\ell : {\bf F}^{\times} \ {\buildrel \sim \over\to}\
K^M_1({\bf F})$; and for $n\geq 2$, $K^M_n({\bf F})$ is presented
as an abelian group by generators $\ell(a_1),\ldots,\ell(a_n)$,
$a_1,\ldots,a_n\in {\bf F}^{\times}$ subject to
\bigskip
$(R_1)_n\hskip1in (a_1,\ldots,a_n) \mapsto \ell(a_1)\cdots\ell(a_n)$
a multilinear function ${\bf F}^{\times}\times \cdots\times
{\bf F}^{\times} \to K^M_n({\bf F})$;
\bigskip
$(R_2)_n \hskip1in \ell(a_1)\cdots\ell(a_n) = 0$ if $a_i + a_{i+1} = 1$ for
some $i < n$.
\bigskip
\underbar{Example} 1.1.  It is customary to express $R_1 \to R_4$ in terms
of the symbol notation:
$$
\{a_1a_2,b\} = \{a_1,b\}\{a_2,b\},
$$
$$
\{a,1-a\} = 1 \ {\rm for}\ a\ne 0,1,
$$
$$
\{a,b\} = \{b,a\}^{-1},
$$
$$
\{a,-a\} = 1.
$$

Furthermore, one can easily verify that
$$
\{a,a\} = \{-1,a\}=\{a,a^{-1}\} = \{a^{-1},a\}.
$$

Now continuing as in [B-T], we introduce $\kappa$-Algebras.
First, a graded ring $\kappa = \coprod_{n\geq 0}
\kappa_n$ is defined by
$$
\kappa = {{\bf Z}[t]\over 2t{\bf Z}[t]} \ = {\bf Z}[\epsilon],
$$
where $\epsilon$ (of degree $1$) $=$ image of $t$. Thus
$\kappa_0 = {\bf Z}$, and for $n\geq 1$, $\kappa_n =
{\bf Z}_2\epsilon^n$, (${\bf Z}_2 = {\bf Z}/2{\bf Z}$). In other
words, $\kappa =$ ring of polynomials in the variable
$\epsilon$ with constant term ${\bf Z}$ and higher degree terms
in ${\bf Z}_2$.
 
\bigskip
{\bf Definition 1.2.} 
A graded $\kappa$-algebra is a graded ring $A = \coprod_{n\geq 0}
A_n$ equipped with a homomorphism $\kappa \to A$ of graded rings,
defined by $\epsilon \mapsto \epsilon_A\in A_1$, such that
$\epsilon_a\in {\rm Center}(A)$. We call $A$ a $\kappa$-Algebra
if further $A_1$ generates $A$ as a $\kappa$-algebra and
$a^2 = \epsilon_aa$ for all $a\in A_1$.
\bigskip
\underbar{Example} 1.3. For a field ${\bf F}$, the map
$\kappa \to K^M_{\ast}({\bf F})$ given by $\epsilon\mapsto \ell(-1)$
gives $K_{\ast}({\bf F})$ the structure of a $\kappa$-Algebra.
Indeed $\ell(-1)$ is central because $K^M_{\ast}({\bf F})$ is
anti-commutative and $2\ell(-1) = 0$; moreover $(R_3^{\prime})
\Rightarrow a^2 = \epsilon_Aa$ (namely $(\ell(a))^2 = \ell(-1)\ell(a)$).
\bigskip
\underbar{Example} 1.4. The free $\kappa$-Algebra on a generator $\Pi$ is
the $\kappa$-Algebra
$$
\kappa(\Pi) = {\kappa[X]\over (X^2 -\epsilon X)},
$$
where $X$ is an indeterminate of degree $1$ with image $\Pi$
modulo $X^2 - \epsilon X$. Evidently, $\kappa(\Pi)$
is a free $\kappa$-modulo with basis $1,\ \Pi$. For any $\kappa$-Algebra
$A$, put $A(\Pi) = A\otimes_{\kappa}\kappa(\Pi) = A \oplus A\Pi$, a
free left $A$-module with basis $1,\ \Pi$. Thus $A(\Pi)_m = A_m\oplus
A_{m-1}\Pi$.
\bigskip
This time we will assume given ${\bf F}$ a field
with a discrete valuation $\nu : {\bf F}^{\times} \to {\bf Z}$,
with corresponding discrete valuation ring ${\cal O} := \{a \in {\bf F} \ |\ 
\nu(a) \geq 0\}$, where we assign
$\nu (0) = \infty$. Let $\pi\in {\cal O}$ generate the unique maximal
ideal $(\pi)$ (i.e. $\nu(\pi) = 1$), and recall that all other non-zero
ideals are of the form $(\pi^m)$, $m \geq 0$. Note that 
${\bf F}^{\times} = {\cal O}^{\times}\pi^{\bf Z}$ (direct product).
Let ${\bf k} = {\bf k}(\nu)$ be the residue field, and $K^M_{\ast}$ Milnor 
$K$-theory. There is a map ($\ell(\pi) = \Pi$):
$$
d_{\pi} : {\bf F}^{\times} \to (K^M_{\ast}{\bf k}(\nu))(\Pi),
\quad d_{\pi}(u\pi^i) = \ell(\overline{u}) + i\Pi,
$$
with  $\overline{u} \in {\bf k}(\nu)$ the corresponding value,
and where $\Pi$ satisfies $\Pi^2 = \ell(-1)\Pi$. This induces:
$$
\partial_{\pi} : K^M_{\ast}{\bf F} \to (K^M_{\ast}{\bf }{\bf k}(\nu))(\Pi).
$$
One next defines maps:
$$
\partial_{\pi}^0,\ \partial_{\nu} : K^M_{\ast}{\bf F} 
\to K^M_{\ast}{\bf k}(\nu),
$$
by
$$
\partial_{\pi}(x) = \partial^0_{\pi}(x) + \partial_{\nu}(x)\Pi,
$$
which can be shown to be independent of the choice of $\pi$
satisfying $\nu (\pi) = 1$ [B-T].
In general $\partial_{\nu} : K^M_m{\bf F} \to K^M_{m-1}{\bf k}(\nu)$. 
Thus for example, suppose we write $\xi = a\pi^i$, then
$
\partial_{\pi}(\xi) = d_{\pi}(\xi) = \ell(\overline{a}) + i\Pi,
$ 
and thus $\partial_{\nu}(\xi) = i$. In general, for $\xi_j = a_j\pi^{k_j}$,
$j = 1,\ldots,m$, we are interested in the product
$$
\partial_{\pi}(\ell(\xi_1)\cdots \ell(\xi_m)) =
(\ell(\overline{a}_1)+k_1\Pi)(\ell(\overline{a}_2)+k_2\Pi)\cdots 
(\ell(\overline{a}_m)+k_m\Pi) = \ell(\overline{a}_1)\cdots 
\ell(\overline{a}_m) + (?)\Pi,
$$
where
$$
\ell(\overline{a}_1)\cdots \ell(\overline{a}_m) =
\partial^0_{\pi}(\ell(\xi_1)\cdots \ell(\xi_m)),\quad 
? = \partial_{\nu}(\ell(\xi_1)\cdots \ell(\xi_m))
$$

A simple calculation gives:
\bigskip
{\bf Proposition 1.5.} 
$$ \partial_{\nu}(\ell(\xi_1)\cdots \ell(\xi_m)) =
\sum_{j=1}^mk_j(-1)^{m-j}\ell(\overline{a}_1)\cdots 
\widehat{\ell(\overline{a}_j)}\cdots \ell(\overline{a}_m) \ 
+\ {\rm intermediate\ terms}\
+ \big(\prod_{j=1}^mk_j\big)\big[\ell(-1)\big]^{m-1}, 
$$
where the ``intermediate terms'' involve the factor $\ell(-1)$.

\vskip.25in
\S {\bf 2. The basic regulator}
\bigskip 
In this section, we define a map
$$
\matrix{r_{\log}&:&H_{\rm Zar}^{k-m}(X,{\cal K}^M_{k,X})&\to&H_{\cal D}^{2k-m}
(X,{\bf R}(k))/(?)\cr
&\cr
&&&&|\wr\cr
&\cr
&&&&{H^{2k-m-1}(X,{\bf R}(k-1))\over
\pi_{k-1}(F^kH^{2k-m-1}(X,{\bf C}))}\big/(?)\cr
&\cr
&&&&\hskip.48in \wr\downarrow\ \matrix{{}_{\rm Hodge}\cr 
{}^{\rm projection}\cr}\cr
&\cr
&&&&\big\{H^{k-1,k-m}(X)\oplus H^{k-m,k-1}(X)\big\}
\bigcap H^{2k-m-1}(X,{\bf R}(k-1))\cr}.
$$
To see how this map is defined, observe that
$$
{H^{2k-m-1}(X,{\bf R}(k-1))\over\pi_{k-1}(F^kH^{2k-m-1}(X,{\bf C}))}
\simeq\big\{H^{k-1,k-m}(X)\oplus\cdots\oplus
H^{k-m,k-1}(X)\big\} \bigcap H^{2k-m-1}(X,{\bf R}(k-1))
$$
$$
\simeq \biggr\{\big\{H^{n-k+m,n-k+1}(X)\oplus\cdots\oplus
H^{n-k+1,n-k+m}(X)\big\}\bigcap H^{2n-2k+m+1}(X,{\bf R}
(n-k+1))\biggr\}^{\vee}
$$
$$
\to \biggr\{\big\{H^{n-k+m,n-k+1}(X)\oplus
H^{n-k+1,n-k+m}(X)\big\}\bigcap H^{2n-2k+m+1}(X,{\bf R}
(n-k+1))\biggr\}^{\vee}.
$$
First, we define a current associated to a basic symbol 
$\{f_1,\ldots,f_m\} \in K^M_m({\bf C}(Z))$,  where codim$_XZ= k-m$.  
Namely, the current defined by:
$$ (2.0)\hskip1in
\sum_{\ell=1}^m\int_Z(-1)^{\ell-1}\log |f_\ell|(d\log |f_1|\wedge\cdots 
\wedge \widehat{d\log |f_\ell|}\wedge\cdots \wedge d\log |f_m|)\wedge w,
$$
where
$$
w\in \big\{E_X^{n-k+m,n-k+1}\oplus
E_X^{n-k+1,n-k+m}\big\}\bigcap \big\{E_X^{2n-2k+m+1}\otimes{\bf R}
(n-k+1))\big\}
$$
is a $C^{\infty}$ differential form on $X$. We must first show
that this current depends only on the symbol 
$\{f_1,\ldots,f_m\} \in K^M_m({\bf C}(Z))$. More specifically,
we prove:
\bigskip
{\bf Proposition 2.1.} Suppose $w$ is both $\partial$ and
$\overline{\partial}$ closed and that $m\geq 2$. 
Then $f_j + f_{j+1} = 1 \ {\rm for\ some}\ j \Rightarrow$
$$
\sum_{\ell=1}^m\int_Z(-1)^{\ell-1}\log |f_\ell|(d\log |f_1|\wedge\cdots 
\wedge \widehat{d\log |f_\ell|}\wedge\cdots \wedge d\log |f_m|)\wedge w = 0.
$$

Proof. Let $F =(f_1,\ldots,f_m) : Z \to \big({\bf P}^1\big)^{\times m}$; 
further let $(t_1,\ldots,t_m)$ be affine coordinates of 
$\big({\bf P}^1\big)^{\times m}$. By a birational modification,
we can assume $Z$ is smooth and that $F$ is a morphism. Now let 
$$
\xi = \sum_{\ell=1}^m(-1)^{\ell-1}\log |t_\ell|(d\log |t_1|\wedge\cdots 
\wedge \widehat{d\log |t_\ell|}\wedge\cdots \wedge d\log |t_m|).
$$
Then
$$
\sum_{\ell=1}^m\int_Z(-1)^{\ell-1}\log |f_\ell|(d\log |f_1|\wedge\cdots 
\wedge \widehat{d\log |f_\ell|}\wedge\cdots \wedge d\log |f_m|)\wedge w = 
\int_Z(F^{\ast}\xi)\wedge w = F_{\ast}(w)(\xi),
$$
where we identify $w$ with its corresponding current. Next, the
assumption $f_j + f_{j+1} = 1 \ {\rm for\ some}\ j \Rightarrow
t_j + t_{j+1} = 1$, and that consequently we can assume given
a morphism $\tilde{F} : Z\to \big({\bf P}^1\big)^{\times m-1}$,
and a corresponding $L^1_{\rm loc}$ form $\tilde{\xi}$ on 
$\big({\bf P}^1\big)^{\times m-1}$, such that 
$$
\sum_{\ell=1}^m\int_Z(-1)^{\ell-1}\log |f_\ell|(d\log |f_1|\wedge\cdots 
\wedge \widehat{d\log |f_\ell|}\wedge\cdots \wedge d\log |f_m|)\wedge w = 
\int_Z(\tilde{F}^{\ast}\tilde{\xi})\wedge w.
$$
But 
$$
\int_Z(\tilde{F}^{\ast}\tilde{\xi})\wedge w = \tilde{F}_{\ast}(w)(\tilde{\xi});
$$
moreover the current $\tilde{F}_{\ast}(w)$ being both $\partial$ and
$\overline{\partial}$ closed implies by $\overline{\partial}$
regularity that $\tilde{F}_{\ast}(w)$ is the current associated
to a holomorphic $m-1$ form on $\big({\bf P}^1\big)^{\times m-1}$
and its conjugate. [Note: The same satement \underbar{cannot} be said of
$F_{\ast}(w)$, for Hodge type reasons!] Since $m\geq 2$, this current must
necessarily be zero, hence the proposition. QED
\bigskip

Now  consider $\eta\in E_X^{n-k+m-1,n-k}$
and $\partial\overline{\partial}\eta + \overline{\partial} 
\partial\overline{\eta}$, and note that $\dim_{\bf R}Z = 2n+2m-2k$.
We need to evaluate, for a given $\{f_1\ldots,f_m\}\otimes Z$,
$$
(2.2)\hskip1in
\sum_{\ell=1}^m\int_Z(-1)^{\ell-1}\log |f_\ell|(d\log |f_1|\wedge\cdots 
\wedge \widehat{d\log |f_\ell|}\wedge\cdots \wedge d\log |f_m|)\wedge
(\overline{\partial}\partial\eta + \partial\overline{\partial}
\overline{\eta}).
$$
First, we can assume by passing to a normalization, that $Z$ is
normal. Further, from (1.1), we can assume that the divisors
for each of the $f_j$, $j=1,\ldots,m$ are different. Let
$\Sigma \subset Z$ be the support of all the divisors $(f_j)$,
$j=1,\ldots,m$, and $\Sigma_{\epsilon}\subset Z$ an $\epsilon$-tube
around $\Sigma$. Now the  integral (2.2) is the same as
$$
\sum_{\ell=1}^m\int_Z(-1)^{\ell-1}\log |f_\ell|(d\log |f_1|\wedge\cdots 
\wedge \widehat{d\log |f_\ell|}\wedge\cdots \wedge d\log |f_m|)\wedge
d(\partial\eta + \overline{\partial}\overline{\eta});
$$
moreover by reason of Hodge type of $\partial\eta + 
\overline{\partial}\overline{\eta}$ and by the taking of residues
along $\Sigma$,
$$
\lim_{\epsilon\to 0}\int_{\Sigma_{\epsilon}}
\log |f_\ell|(d\log |f_1|\wedge\cdots 
\wedge \widehat{d\log |f_\ell|}\wedge\cdots \wedge d\log |f_m|)\wedge
(\partial\eta + \overline{\partial}\overline{\eta}) = 0.
$$
Thus by Stokes' Theorem, the calculation of (2.2) amounts to
calculating
$$
\int_Z(d\log |f_1|\wedge\cdots \wedge d\log |f_m|)\wedge
(\partial\eta + \overline{\partial}\overline{\eta}).
$$
We can further reduce this to a calculation of the type
$$
\int_Z({d\overline{f}_1\over \overline{f}_1}\wedge\cdots 
\wedge {d\overline{f}\over \overline{f}_m})\wedge\partial\eta 
\ + \ \int_Z({df_1\over f_1}\wedge\cdots\wedge {df_m\over f_m})
\wedge\overline{\partial}\overline{\eta}.
$$
Again, by Hodge type considerations, this is the same as
$$
\int_Z({d\overline{f}_1\over \overline{f}_1}\wedge\cdots 
\wedge {d\overline{f}\over \overline{f}_m})\wedge d\eta 
\ + \ \int_Z({df_1\over f_1}\wedge\cdots\wedge {df_m\over f_m})
\wedge d\overline{\eta},
$$
and by Stokes' Theorem, this amounts to calculating
$$
\lim_{\epsilon\to 0} \biggl(
\int_{\Sigma_{\epsilon}}({d\overline{f}_1\over \overline{f}_1}\wedge\cdots 
\wedge {d\overline{f}\over \overline{f}_m})\wedge \eta 
\ + \ \int_{\Sigma_{\epsilon}}({df_1\over f_1}\wedge\cdots
\wedge {df_m\over f_m})\wedge \overline{\eta}\biggr).
$$
By a residue calculation, this amounts to the same thing as
$$
= \sum_{D\subset \Sigma}\biggl(\sum_{\ell=1}^m(-1)^{m-\ell}\nu_{D}(f_{\ell})
\int_D(d\log |f_1|\wedge\cdots 
\wedge \widehat{d\log |f_\ell|}\wedge\cdots \wedge d\log |f_m|)\wedge
(\eta + \overline{\eta})\biggr).
$$
This expression is strikingly familiar to the first term
of  (1.5). Indeed, if in (1.5), we set ${\bf F} = {\bf C}(D)$,
and replace $\ell(\ )$ by $d\log|(\ )|$, then $d\log|(-1)| = 0$
and integration leads precisely to the above calculation. In particular,
if a class in $H^{k-m}_{\rm Zar}(X,\overline{\cal K}^M_{k,X})$
is represented by $\zeta$, then the corresponding current $T_\zeta$ associated
to $\zeta$, induced by the formula in (2.0), is 
$\partial\overline{\partial}$-closed. Next, suppose 
$w\in E_X^{n-k+m,n-k+1} \oplus E_X^{n-k+1,n-k+m}$ is a real $\partial$
and $\overline{\partial}$ closed form.  For a given
$\{f_1,\ldots,f_{m+1}\}\otimes Z$, we will again assume that
$f_j \in {\bf C}(Z)^{\times}$, the divisors of $(f_j)$, 
$j=1,\ldots,m+1$ are all different, and that we have a 
morphism $F = (f_1,\ldots,f_{m+1}) : Z\to \big({\bf P}^1\big)^{\times m +1}$, 
where codim$_XZ = k-m-1$, and $Z$ is presumed smooth.
Thus we can write
$\{f_1,\ldots,f_{m+1}\}\otimes Z = F^{\ast}\{t_1,\ldots,t_{m+1}\}
\otimes ({\bf P}^1)^{\times m+1}$, where $(t_1,\ldots,t_{m+1})$
are affine coordinates of $\big({\bf P}^1\big)^{\times m+1}$. For
a given irreducible divisor $D\subset \big({\bf P}^1\big)^{\times m+1}$,
we consider the elements
$$
\xi_D := \sum_{j=1}^{m+1}(-1)^{m+1-j}\nu_D(t_j)
\{t_1,\ldots,\hat{t}_j,\ldots,t_{m+1}\}_D,
\quad \xi := \sum_D\xi_D.
$$
$$ 
\eta_{\xi_D} := 
\sum_{j=1}^{m+1}(-1)^{m+1-j}\nu_D(t_j)d\log|t_1|\wedge\cdots
\wedge \widehat{d\log|t_j|}\wedge\cdots\wedge d\log|t_{m+1}|.
$$
If we can show that
$$(2.3)\hskip2in
\sum_D\int_DF^{\ast}(\eta_{\xi_D})\wedge w = 0,
\hskip2in
$$
then by (1.5) and functoriality, we arrive at a well-defined map
$$
H^{k-m}_{\rm Zar}(X,\overline{\cal K}^M_{k,X}) \to 
\biggl\{\big\{H^{n-k+m,n-k+1}(X)\oplus H^{n-k+1,n-k+m}\big\}
\bigcap H^{2n-2k+m+1}(X,{\bf R}(n-k+1))\biggr\}^{\vee},
$$
hence by duality a map
$$
H^{k-m}_{\rm Zar}(X,\overline{\cal K}^M_{k,X}) \to
\big\{H^{k-1,k-m}(X)\oplus H^{k-m,k-1}(X)\big\} 
\bigcap H^{2k-m-1}(X,{\bf R}(k-1)).
$$
The only irreducible codimension one $D$'s in 
$\big({\bf P}^1\big)^{\times m+1}$ contributing to non-trivial
$\xi_D$ are the following. Set $D_{0,\ell} = 
{\bf P}^1\times\cdots\times \{0\}\times \cdots
\times {\bf P}^1$, $D_{\infty,\ell} = {\bf P}^1\times\cdots\times 
\{\infty\}\times \cdots\times {\bf P}^1$. Then
$$
\xi_{D_{0,\ell}} = (-1)^{m+1-\ell}\{t_1,\ldots,
\hat{t}_{\ell},\ldots,t_{m+1}\}_{D_{0,\ell}}
$$
$$
\xi_{D_{\infty,\ell}} = (-1)^{m-\ell}\{t_1,\ldots,
\hat{t}_{\ell},\ldots,t_{m+1}\}_{D_{\infty,\ell}}.
$$
Note that $D_{0,\ell}\simeq ({\bf P}^1)^{\times m} \simeq D_{\infty,\ell}$,
and that $D_{0,\ell} \sim_{\hom} D_{\infty,\ell}$ in 
$({\bf P}^1)^{\times m+1}$. Now $\xi$ defines a 
$\partial\overline{\partial}$-closed current $T_\xi$ 
of the form $\sum_D\int_DG_D$, where $G_D$ is a $\log$ form 
that is $\partial\overline{\partial}$ closed. 
By the Poincar\'e-Lelong Theorem, we can write the current of
integration over $D$, namely $\delta_D := \int_D = \partial
\overline{\partial}T_{\psi_D} + \Phi_D$, where $\psi_D$ is $\log$ type,
and $\Phi_D$ is a $\partial$ and $\overline{\partial}$ closed
$C^{\infty}$ form on $({\bf P}^1)^{\times m+1}$, namely the
first Chern form. Note that for a $C^{\infty}$ bump function 
$\varphi_{\epsilon}$ on $({\bf P}^1)^{m+1}$, zero on an 
$\epsilon$-neighbourhood $\Sigma_{<\epsilon}$ of the
various $D$'s in question,
$$
\partial\overline{\partial}T_{\psi_D}(G_D) = \lim_{\epsilon \to 0}
\partial\overline{\partial}T_{\varphi_{\epsilon}\cdot\psi_D}(G_D)
= \lim_{\epsilon \to 0}T_{\varphi_{\epsilon}\cdot\psi_D}
(\partial\overline{\partial}G_D) = 0.
$$
Thus $T_\xi = \sum_D\int_{({\bf P}^1)^{\times m+1}}G_D\wedge \Phi_D$.
If we work with $D = D_{0,\ell}$ say, then 
$$
\Phi_D = {1\over 2\pi\sqrt{-1}}{dt_{\ell}\wedge d\overline{t}_{\ell}
\over (1 + |t_{\ell}|^2)^2}.
$$
This is the same form for $D_{\infty,\ell}$ since the Chern
form depends only on the bundle. Also, as forms on $({\bf P}^1)^{\times m+1}$,
$G_{D_{0,\ell}} + G_{D_{\infty,\ell}} = 0$. For each $D$,
we are essentually dealing with the calculation
$$
G_D\wedge \Phi_D = \sum_r(-1)^*\log |t_r|{1\over 2\pi\sqrt{-1}}
{dt_{\ell}\wedge d\overline{t}_{\ell}\over (1 + |t_{\ell}|^2)^2}
\wedge \biggl(\bigwedge_{j=1, j\ne \ell, r}^{m+1}d\log |t_j|\biggr).
$$
Moreover, for any current $T$ satisfying 
$\partial\overline{\partial} T = F_{\ast}w - T_{\eta}$ for
some $C^{\infty}$ closed form $\eta$, and if $\varphi_{\epsilon}$ is
a $C^{\infty}$ bump function vanishing on $\Sigma_{<\epsilon}$, then
$$
(2.4)\hskip1in 
\partial\overline{\partial}T(G_D\wedge \Phi_D) = \lim_{\epsilon\to 0}
\partial\overline{\partial}T(\varphi_{\epsilon}G_D\wedge \Phi_D) 
= \lim_{\epsilon\to 0}T(\partial\overline{\partial}
(\varphi_{\epsilon}\cdot\wedge G_D\wedge \Phi_D)),
$$ 
and by symmetry considerations
together with $G_{D_{0,\ell}} + G_{D_{\infty,\ell}} = 0$, 
it follows that $\sum_D\partial\overline{\partial}T
\big(G_D\wedge \Phi_D\big) = 0$.
It follows that $F_{\ast}w$ in the formula
$$
\sum_D\int_DF^{\ast}(\eta_{\xi_D})\wedge w =
\sum_D\int_{{{\bf P}^1}^{m+1}}G_D\wedge \Phi_D\wedge F_{\ast}(w)
$$
can be replaced by any $C^{\infty}$ closed form 
representative $\eta$, where 
$$
\{\eta\} \in H^{m+1,0}(\big({\bf P}^1\big)^{\times m +1})\oplus
H^{0,m+1}(\big({\bf P}^1\big)^{\times m +1}) = 0.
$$  
But such an $\eta$
can thus be chosen to be zero. QED
\vskip.25in
\S {\bf 3. Comparison to the Beilinson regulator}
\bigskip
According to [Lew3], the definition of real Deligne cohomology
of a smooth quasi-projective variety $U$ with good
compactification $\overline{U}$ (with $E := \overline{U}
\backslash U$) can be taken to be given by:
$$
H^i_{\cal D}(U,{\bf R}(p)) \simeq H^i\big(\tilde{{\bf R}(p)}_{\cal D}
(\overline{U}) :=
{\rm Cone}\{F^p\Omega^{\bullet}_{\overline{U}^{\infty}}<E>(\overline{U})  \
{\buildrel -\pi_{p-1}\over\longrightarrow} \
{\cal E}^{\bullet}_{\overline{U},\bf R}<E>(p-1)(\overline{U})\}[-1]\big),$$
where ${\cal E}^{\bullet}_{\overline{U},\bf R}<E>(\overline{U})$
is a corresponding real logarithmic complex [Bu].
The product structure on $\tilde{{\bf R}(p)}_{\cal D}$,
viz: 
$$\tilde{{\bf R}(p)}_{\cal D} \times
\tilde{{\bf R}(q)}_{\cal D} \ {\buildrel \cup \over \longrightarrow}
\ \tilde{{\bf R}(p+q)}_{\cal D}$$
can be arrived at from the table in [E-V], and is given below:
\bigskip
$$
\matrix{&\vert&f_q&\vert&s_q\cr
---&&---&&---\cr
f_p&\vert&f_p\wedge f_q&\vert&(-1)^{\deg f_p}\pi_p(f_p)\wedge s_q\cr
---&&---&&---\cr
s_p&\vert&s_p\wedge\pi_q(f_q)&\vert&0\cr}
$$
 
This defines also defines  multiplication on 
the corresponding description of $H^i_{\cal D}(X,{\bf R}(p))$
via the definition above. For instance, $H_{\rm Zar}^0(U,
\overline{K}^M_{1,U}) = H^0(U,{\cal O}_U^{\times})$, and the
corresponding Beilinson regulator map $r_B : H^0(U,{\cal O}_U^{\times})
\to H^1_{\cal D}(U,{\bf R}(1))$ is represented by
$$
f\in H^0(U,{\cal O}_U^{\times}); \quad r_B(f) = \{({df\over f},\log |f|)\}
\in H^1_{\cal D}(U,{\bf R}(1)). 
$$
Now let $f_1,\ldots,f_m\in {\bf C}(U)^{\times}$ be given.
Then we have:
$$
\{({df_1\over f_1},\log |f_1|)\}\bigcup \cdots \bigcup 
\{({df_m\over f_m},\log |f_m|)\}
$$
$$
= \biggl\{\biggl({df_1\over f_1}\wedge\cdots\wedge {df_m\over f_m},
\xi(f_1,\ldots,f_m)\biggr)\biggr\}\in H^m_{\cal D}(U,{\bf R}(m)),
$$
where for example
$$
\xi(f_1) = \log|f_1|;\quad \xi(f_1,f_2) = \log|f_1|\pi_1\big(
{df_2\over f_2}\big) - \log|f_2|\pi_1\big({df_1\over f_1}\big);
$$
$$
\xi(f_1,f_2,f_3) = \log|f_1|\pi_1\big({df_2\over f_2}\big)\wedge
\pi_1\big({df_3\over f3}\big) - \log|f_2|\pi_1\big({df_1\over f_1}\big)
\wedge\pi_1\big({df_3\over f_3}\big) + \log|f_3|
\pi_2\big({df_1\over f_1}\wedge{df_2\over f_2}\big);
$$
$$
\xi(f_1,f_2,f_3,f_4) = \log|f_1|\pi_1\big({df_2\over f_2}\big)\wedge
\pi_1\big({df_3\over f_3}\big)\wedge \pi_1\big({df_4\over f_4}\big)
- \log|f_2|\pi_1\big({df_1\over f_1}\big)\wedge\pi_1\big({df_3\over f_3}\big)
\wedge \pi_1\big({df_4\over f_4}\big)
$$
$$ 
+ \log|f_3|\pi_2\big({df_1\over f_1}\wedge{df_2\over f_2}\big)\wedge
\pi_1\big({df_4\over f_4}\big) - \log|f_4|\pi_3\big(
\big({df_1\over f_1}\wedge{df_2\over f_2}\big)\wedge \pi_1
\big({df_3\over f_3}\big),
$$
 
and so on $\ldots$
\bigskip
We want to apply this to the following setting. Fix
an irreducible subvariety $Z \subset X$ of codimension $k-m$,
and $f_1,\ldots,f_m\in {\bf C}(Z)^{\times}$.
Let $D \subset Z$ be the divisor supporting the zeros
and poles of the $f_j$'s, $j=1,\ldots,m$ (and $Z_{\rm sing})$. 
Let $U_X = X\backslash D$
and $U_Z = Z \backslash D$. We refer to the following diagram.
$$
\matrix{H_{\rm Zar}^{k-m}(X,\overline{\cal K}^M_{k,X})&\to&
\{H^{k-1,k-m}(X)\oplus H^{k-m,k-1}(X)\}\cap
H^{2k-m-1}(X,{\bf R}(k-1))\cr
&\cr
\swarrow\hskip.2in &&\uparrow\cr
&\cr
E_2^{k-m,-k}(k)\Rightarrow CH^k(X,m)\hskip.2in
&\to&H_{\cal D}^{2k-m}(X,{\bf R}(k))\cr
&\cr
\hskip.2in \swarrow&&\downarrow\cr
&\cr
CH^k(U_X,m)&\to&H_{\cal D}^{2k-m}(U_X,{\bf R}(k))\cr
&\cr
\uparrow&&\quad\uparrow i_{\ast}\cr
&\cr
CH^m(U_Z,m)&\to&H_{\cal D}^{m}(U_Z,{\bf R}(m))\cr
&\cr
\cap\uparrow\quad&&\quad\uparrow\cup\cr
&\cr
CH^1(U_Z,1)^{\otimes m}&\to&H_{\cal D}^{1}(U_Z,{\bf R}(1))^{\otimes m}\cr
\cr}
$$

We prove:
\bigskip
{\bf Proposition 3.0.} The Beilinson regulator\footnote{$^{\dagger}$}
{Wherever $r_B$ is defined. Let us refer back to the spectral
sequence in \S 0. Note that in 
general for $\ell \geq 1$ and any $m\geq 0$,
$E_2^{k-m+\ell,-k-\ell} = H_{\rm Zar}^{k-m+\ell}(X,{\cal CH}^k(k
+\ell))$. From the Gersten resolution, this involves terms of 
the form $CH^{m-\ell}({\bf C}(Z),m)$, where codim$_XZ = k-m+\ell$, i.e.
$\dim Z = n-k+m-\ell$, which cannot support a form $w$ in
$E_X^{n-k+m,n-k+1}\oplus E_X^{n-k+1,n-k+m}$, since $\ell \geq 1$.
Thus, so long as one can represent a class $\xi\in E_2^{k-m,-k}
= H^{k-m}_{\rm Zar}(X,\overline{\cal K}^M_{k,X})$ by
a lifting $\tilde{\xi}\in CH^k(X,m)$ (for example if
$\xi$ lives forever in this spectral sequence), then $r_B(\xi)$
is defined and given by $r_B(\xi) = r_{\cal D}(\tilde{\xi})$.}
$$
r_B : H_{\rm Zar}^{k-m}(X,\overline{\cal K}^M_{k,X})
\to \{H^{k-1,k-m}(X)\oplus H^{k-m,k-1}(X)\} \bigcap
H^{2k-m-1}(X,{\bf R}(k-1))
$$
is induced by 
$
w \in \{H^{n-k+m,n-k+1}(X)\oplus H^{n-k+1,n-k+m}(X)\} \bigcap
H^{2k-m-1}(X,{\bf R}(n-k+1))
$
$$ \mapsto
{1\over (2\pi\sqrt{-1})^{n-k+m}}\sum_{\ell=1}^m\int_Z
\xi(f_1,\ldots,f_m)\wedge w.
$$

Proof. Let $\overline{U}_X$ be a good
compactification of $U_X$. In particular, we have a morphism
$\sigma : \overline{U}_X \to X$, where $\sigma^{-1}(D)
= \overline{D}$ is a NCD in $\overline{U}_X$. Then $\big(
\overline{U}_X,\overline{D}\big)$ can be used to compute the
Deligne cohomology $H_{\cal D}^{2k-m}(U_X,{\bf R}(k))$.
Let $\zeta = \sum_{\alpha}\big\{f_{1,\alpha},\ldots,
f_{m,\alpha}\big\}\otimes Z_{\alpha}\in H_{\rm Zar}^{k-m}(X,
\overline{\cal K}^M_{k,X})$ be given, with corresponding image
in $r_B(\zeta) \in H_{\cal D}^{2k-m}(X,{\bf R}(k))$. Likewise,
we have a corresponding $i_{\ast}({df\over f},\xi) := 
i_{\ast}\big(\sum_{\alpha}({df_{1,\alpha}\over f_{1,\alpha}}
\wedge {df_{m,\alpha}\over f_{m,\alpha}},\xi(f_{1,\alpha},
\ldots,f_{m,\alpha}))\big) \in H_{\cal D}^{2k-m}(U_X,{\bf R}(k))$. 
Then working over $U_X$, $r_B(\zeta) - i_{\ast}({df\over f},\xi)$ 
is a coboundary in $H_{\cal D}^{2k-m}(U_X,{\bf R}(k))$. 
Working with the second factor in the cone description 
of $H_{\cal D}^{2k-m}(U_X,{\bf R}(k))$, and using $\sigma^{\ast}
(E_X^{n-k+m,n-k+1}\oplus E_X^{n-k+1,n-k+m}) \subset 
\Sigma_{\overline{D}}E^{\bullet}_{\overline{U}_X} :=$ forms
which pull back to zero on $\overline{D}$, it follows from 
the techniques used in [Lew1] that 
$\big(r_B(\zeta) - i_{\ast}({df\over f},\xi)\big)
(w) = 0$\footnote{$^{\dagger\dagger}$}{Strickly speaking,
$r_B(\zeta) - i_{\ast}({df\over f},\xi)$ is a coboundary 
\underbar{current}.
But that current can be chosen in the dual space of
$\Sigma_{\overline{D}}E^{\bullet}_{\overline{U}_X}$
because $\overline{U}_X$ is smooth and $\overline{D}
\subset \overline{U}_X$ is a NCD. See [Bu](Remark,
page 562).}. Part
(iii) of the Main Theorem is easy and will be left to the reader.

\vfill\eject
{\bf References}
\bigskip
\item{[B-T]} H. Bass and J. Tate, The Milnor ring of a global
field, in Algebraic $K$-theory II, Lecture Notes in Math. {\bf 342},
Springer-Verlag, 1972, 349--446.

\item{[Bei]} A. Beilinson, Higher regulators and values of
$L$-functions, J. Soviet math. {\bf 30}, 1985, 2036--2070.

\item{[Blo]} S. Bloch, Algebraic cycles and higher $K$-theory,
Adv. Math. {\bf 61}, 1986, 267--304.

\item{[Bu]} J. Burgos, Green forms and their product, Duke
Math. J. {\bf 75}, No. {\bf 3}, 1994, 529--574.

\item{[EV]} H. Esnault and E. Viehweg, Deligne-Beilinson cohomology,
in Beilinson's Conjectures on Special Values of $L$-Functions, (Rapoport,
Schappacher, Schneider, eds.), Perspect. Math. {\bf 4}, Academic Press, 
1988, 43--91.

\item{[G-L]} B. Gordon and J. Lewis, Indecomposable higher Chow cycles
on products of elliptic curves, J. Algebraic Geometry {\bf 8}, 
1999, 543--567.

\item{[GH]} P.
Griffiths and J. Harris, Principles of Algebraic Geometry, John Wiley
\& Sons, New York, 1978.

\item{[Ja]} U. Jannsen, Deligne cohomology, Hodge-${\cal D}$-conjecture,
and motives, in Beilinson's Conjectures on Special Values of $L$-Functions, 
(Rapoport, Schappacher, Schneider, eds.), Perspect. Math. {\bf 4}, 
Academic Press, 1988, 305--372.

\item{[Ka]} K. Kato, Milnor $K$-theory and the Chow group of zero
cycles, Contemporary Mathematics, Vol. {\bf 55}, Part I, 1986, 241--253.

\item{[Lev]} M. Levine, Localization on singular varieties, 
Invent. Math. {\bf 31}, 1988, 423--464.

\item{[Lew]} Lewis, J.D.: Higher Chow groups and the Hodge-$\cal D-$conjecture,
Duke Math. J. {\bf 85} (1996), 183--207.

\item{} \vrule height0pt width70pt, 
A note on indecomposable motivic cohomology classes,
J. reine angew. Math. {\bf 485} (1997), 161--172.

\item{} \vrule height0pt width70pt,
A duality pairing between cohomology and higher Chow groups,
J. reine angew. Math. {\bf 504} (1998), 177--193.

\item{[MS]} S. M\"uller-Stach, Constructing indecomposable motivic 
cohomology classes on algebraic surfaces, J. Algebraic Geometry {\bf 6},
1997, 513--543. 

\item{} \vrule height0pt width70pt, Algebraic cycle complexes, in Proceedings
of the NATO Advanced Study Institute on the Arithmetic and
Geometry of Algebraic Cycles Vol. {\bf 548}, 
(Lewis, Yui, Gordon, M\"uller-Stach, 
S. Saito, eds.), Kluwer Academic Publishers, Dordrecht, The 
Netherlands, (2000), 285--305.

\item{[Sou]} C. Soul\'e, Lectures on Arakelov Geometry, Cambridge Studies
in Advanced Mathematics {\bf 33}, 1992.

\bye

\vfill\eject
${\cal CH}^k(r) =$ sheaf associated to the presheaf
$U\subset X$ open $\mapsto CH^k(U,r)$. Fix an integer $k \geq 0$.
There is a spectral sequence, $E_2^{p,q} := H^p_{\rm Zar}(X,{\cal CH}^k(-q))
\Rightarrow CH^k(X,-p-q)$. We now set $p+q=-m$, or $-q = p+m$.
Note that $p\geq 0,\ q\leq 0$. So $E_2^{p,-p-m} =
H^p(X,{\cal CH}^k(p+m))\Rightarrow CH^k(X,m)$. Further,
${\cal CH}^k(k) \simeq \overline{\cal K}^M_{k,X}$, and one expects
in characteristic zero that $\overline{\cal K}^M_{k,X} =
{\cal K}^M_{k,X}$. Now if we set $k=p+m$, then
$E_2^{k-m,-k} = H^{k-m}_{\rm Zar}(X,\overline{\cal K}^M_{k,X})$.
Bloch proves the existence of a Gersten resolution
$$
0\to {\cal CH}^k(r)\to \bigoplus_{x\in X^0}i_xCH^k({\rm Sp}(k(x)),r)\to
\bigoplus_{x\in X^1}i_xCH^{k-1}({\rm Sp}(k(x)),r-1)\to\cdots
$$
$$
\to \bigoplus_{x\in X^r}i_xCH^{k-r}({\rm Sp}(k(x)),0)\to 0
$$

\underbar{Facts}: Let ${\bf F}$ be a field, and identify
with its prime spectrum.
\bigskip
(i) $CH^\nu({\bf F},\mu) = 0$ if $\nu > \mu$
\bigskip
(ii) $CH^\nu({\bf F},\mu) = 0$ if $\nu = 0$ \underbar{and} $\mu \geq 1$;
or  if $\nu = 1$ \underbar{and} $\mu \geq 2$.
\bigskip
But $E_2^{k-m+\ell,-k-\ell} = H_{\rm Zar}^{k-m+\ell}(X,{\cal CH}^k(k
+\ell))$ and $k-(k-m+\ell) = m-\ell$, $k+\ell - (k-m+\ell) = m$.
Thus via the Gersten resolution, we are focussing on $CH^{m-\ell}({\bf F},
m)$ which is zero if $m\leq 2$ (with $\ell \geq 1$) by fact (ii) above.
\bigskip
Next, if $r\geq 2$, we consider $E_2^{k-m+r,-k-r+1} = H^{k-m+r}_{\rm Zar}
(X,{\cal CH}^k(k+r-1))$. Again, $k - (k-m+r) = m-r,\
k+r-1-(k-m+r) = m-1$. We are looking at $CH^{m-r}({\bf F},m-1)$ which
is zero for $m\leq 2$ ($r\geq 2$). Therefore there is an injection
$\lambda : H_{\rm Zar}^{k-m}(X,\overline{\cal K}^M_{k,X})
\hookrightarrow CH^k(X,m)$. Finally, I claim that $\lambda$ is
an isomorphism for $m\leq 2$. We assume $\ell \geq 1$ and look
at $E_2^{k-m-\ell,-k+\ell} = H_{\rm Zar}^{k-m-\ell}(X,
{\cal CH}^k(k-\ell))$, but this is zero, using
fact (i) above and the Gersten resolution. 
\bigskip
Note further that in general for $\ell \geq 1$ and any $m\geq 0$,
$E_2^{k-m+\ell,-k-\ell} = H_{\rm Zar}^{k-m+\ell}(X,{\cal CH}^k(k
+\ell))$ and $k-(k-m+\ell) = m-\ell$, $k+\ell - (k-m+\ell) = m$.
Thus via the Gersten resolution, we are focussing on 
$CH^{m-\ell}({\bf C}(Z),m)$, where codim$_XZ = k-m+\ell$, i.e.
$\dim Z = n-k+m-\ell$, whch cannot support a form $w$ in
$E_X^{n-k+m,n-k+1}\oplus E_X^{n-k+1,n-k+m}$, since $\ell \geq 1$.
Thus, so long as one can represent a class $\xi\in E_2^{k-m,-k}$ by
a by a lifting $\tilde{\xi}\in CH^k(X,m)$ (for example if
$\xi$ lives forever in this spectral sequence), then $r_B(\xi)$
is defined and given by $r_B(\xi) = r_{\cal D}(\tilde{\xi})$.

\vfill\eject

$\dim Z = n-k+2$, $\xi,\ w\in E_X^{n-k+2,n-k+1}$,
$f, g\in {\bf C}(Z)^{\times}$. Introduce:
$$
r_{\log}(f,g)(w) = \int_Z\big[\log|g|\big({df\over f} + {d\overline{f}\over
\overline{f}}\big) - \log|f|\big({dg\over g} +
{d\overline{g}\over\overline{g}}\big)\big]\wedge(w\oplus \overline{w})
$$
$$
= \int_Z\big[\log|g|{d\overline{f}\over\overline{f}} -
\log|f|{d\overline{g}\over\overline{g}}\big]\wedge w\ +\
\int_Z\big[\log|g|{df\over f} - \log|f|{dg\over g}\big]\wedge
\overline{w}
$$
$$
R_D(f,g)(w) = \int_Z\big[\log|g|\big({df\over f} - {d\overline{f}\over
\overline{f}}\big) - \log|f|\big({dg\over g} -
{d\overline{g}\over\overline{g}}\big)\big]\wedge(\xi\ominus \overline{\xi})
$$
$$
= -\int_Z\big[\log|g|{d\overline{f}\over\overline{f}} -
\log|f|{d\overline{g}\over\overline{g}}\big]\wedge \xi\ -\
\int_Z\big[\log|g|{df\over f} - \log|f|{dg\over g}\big]\wedge
\overline{\xi}
$$
Both integrals are real valued. They compare by the formula:
$$
r_{\log}(f,g)(w\oplus\overline{w}) = -R_D(f,g)(w\ominus \overline{w})
$$
We now have:
$$
r_D(f,g)\big((2\pi\sqrt{-1})^{n-k+1}(w\oplus\overline{w})\big)
:= {1\over (2\pi\sqrt{-1})^{n-k}}R_D(f,g)\big((2\pi\sqrt{-1})^{n-k+1}
(w\oplus\overline{w})\big)
$$
$$
= R_D(f,g)(2\pi\sqrt{-1}[w\oplus\overline{w}]) = R_D\big((2\pi\sqrt{-1}w)
\ominus (\overline{2\pi\sqrt{-1}w})\big)
$$
$$
= r_{\log}(f,g)\big((2\pi\sqrt{-1}w)\oplus (\overline{2\pi\sqrt{-1}w})\big)
$$
$$
= {1\over (2\pi\sqrt{-1})^{n-k+1}}r_{\log}(f,g)\big((2\pi\sqrt{-1})^{n-k+1}
[(2\pi\sqrt{-1}w)\oplus (\overline{2\pi\sqrt{-1}w})]\big)
$$
$$
:= r_{\log}(f,g)\big((2\pi\sqrt{-1})^{n-k+1}[(2\pi\sqrt{-1}w)
\oplus (\overline{2\pi\sqrt{-1}w})]\big)
$$
Thus we have compared these regulators. 
{\bf Proposition 2.0.} Suppose $w$ is both $d$-closed. 
Then $f_j + f_{j+1} = 1 \Rightarrow$
$$
\sum_{\ell=1}^m\int_Z(-1)^{\ell-1}\log |f_\ell|(d\log |f_1|\wedge\cdots 
\wedge \widehat{d\log |f_\ell|}\wedge\cdots \wedge d\log |f_m|)\wedge w = 0.
$$

Proof. Let $D \subset X$ be the divisor supporting the zero
and pole of the $f_j$'s, $j=1,\ldots,m$. Let $U_X = X\backslash D$
and $U_Z = Z \cap U_X$. The following diagram will be referred to
here, and for the remainder of the paper.
$$
\matrix{&&H_{\rm Zar}^{k-m}(X,\overline{\cal K}^M_{k,X})\cr
&\cr
&&\downarrow&\searrow\cr
&\cr
&&CH^k(X,m)&\to&H_{\cal D}^{2k-m}(X,{\bf R}(k))\cr
&\cr
&&\downarrow&&\downarrow\cr
&\cr
&\nearrow&CH^k(U_X,m)&\to&H_{\cal D}^{2k-m}(U_X,{\bf R}(k))\cr
&\cr
&&\uparrow&&\uparrow\cr
&\cr
CH^m(Z,m)&\to&CH^m(U_Z,m)&\to&H_{\cal D}^{m}(U_Z,{\bf R}(m))\cr
&\cr
&&\cap\uparrow\quad&&\quad\uparrow\cup\cr
&\cr
&&CH^1(U_Z,1)^{\otimes m}&\to&H_{\cal D}^{1}(U_Z,{\bf R}(1))^{\otimes m}\cr
\cr}
$$
By construction, the class $\xi \in H_{\rm Zar}^{k-m}
(X,\overline{\cal K}^M_{k,X})$ arising from the symbol
$\{f_1,\ldots,f_m\}\otimes Z$ (with $f_j + f_{j+1} = 1$
for some $j$) is both  supported on $Z$ and zero. By a diagram
chase and by the work of Burgos (and my recent Crelle paper),
the current defined in (2.0) is the coboundary current
defined by the differential of a form with log singularities, modulo
$F^m$, which acts as zero on cohomology.
\bye